# NONPARAMETRIC METHODS FOR INFERENCE IN THE PRESENCE OF INSTRUMENTAL VARIABLES

By Peter Hall and Joel L. Horowitz[1]

*Australian National University and Northwestern University*

We suggest two nonparametric approaches, based on kernel methods and orthogonal series to estimating regression functions in the presence of instrumental variables. For the first time in this class of problems, we derive optimal convergence rates, and show that they are attained by particular estimators. In the presence of instrumental variables the relation that identifies the regression function also defines an ill-posed inverse problem, the "difficulty" of which depends on eigenvalues of a certain integral operator which is determined by the joint density of endogenous and instrumental variables. We delineate the role played by problem difficulty in determining both the optimal convergence rate and the appropriate choice of smoothing parameter.

**1. Introduction.** Data $(X_i, Y_i)$ are observed, the pairs being generated by the model

$$Y_i = g(X_i) + U_i, \qquad (1.1)$$

where $g$ is a function which we wish to estimate and the $U_i$'s denote disturbances. The $U_i$'s are correlated with the explanatory variables $X_i$ and, in particular, $E(U_i|X_i)$ does not vanish. For example, this may occur if a third variable causes both $X_i$ and $Y_i$, but is not included in the model.

This circumstance arises frequently in economics. To illustrate, suppose that $Y_i$ denotes the hourly wage of individual $i$, and that $X_i$ includes the individual's level of education, among other variables. The "error" $U_i$ would generally include personal characteristics, such as "ability," which influence the individual's wage but are not observed by the analyst. If high-ability

Received March 2003; revised December 2004.
[1]Supported in part by NSF Grants SES-99-10925 and SES-03-52675.

individuals tend to choose high levels of education, then education is correlated with ability, thereby causing $U_i$ to be correlated with at least some components of $X_i$.

Suppose, however, that for each $i$ we have available another observed data value, $W_i$, say (an instrumental variable), for which

$$(1.2) \qquad E(U_i|W_i) = 0$$

and there is a "sufficiently strong" relationship between $X_i$ and $W_i$. Then there is an opportunity for estimating $g$ from the data $(X_i, W_i, Y_i)$.

The formal definition of "sufficiently strong" will depend on the nature of the problem. In a parametric setting, for example, where $g(X_i) = X_i\beta$, $X_i$ is an $m \times k$ matrix and $\beta$ is a $k \times 1$ vector, "sufficiently strong" means simply that the matrix of correlations between $X$ and $W$ is of full rank; this is sometimes expressed as "$X$ and $W$ are fully correlated." In a nonparametric setting the definition of "sufficiently strong" is given by, for example, condition (2.1) below.

Estimation of $g$ is difficult because, as explained in Section 2, the relation that identifies $g$ is a Fredholm equation of the first kind,

$$(1.3) \qquad Tg = \phi,$$

say, which leads to an ill-posed inverse problem [9, 14]. We use a ridge-type regularization method to achieve boundedness of the relevant inverse integral operator, and develop both kernel and series estimators of $g$. The resulting estimators have optimal $L_2$ rates of convergence.

Closely related inverse problems, where the context is rendered relatively abstract in order to facilitate solution, include those studied by Donoho [4], Johnstone [8] and Cavalier, Golubev, Picard and Tsybakov [2]. That work addresses the white-noise model, rather than the more explicitly realistic discrete-data setting of (1.1). In such treatments the operator $T$ is generally assumed known, whereas in the case of instrumental-variables problems it usually must be estimated from data. Nevertheless, the optimal convergence rates obtained in the above earlier work are identical to our own. Indeed, the mean integrated squared error rates we obtain are the same as those in an "ordinary" inverse problem, where $T$ is known and equal to $T_1'T_1$, and $T_1$ is the nonstochastic transformation of the actual inverse model. Efromovich and Koltchinskii [5] treated a white-noise model in a setting where $T$, at (1.3), must be estimated, and also obtained optimal rates.

Research on this type of problem in econometrics is mostly very recent. Blundell and Powell [1] and Florens [6] discussed the relationship between (1.1) and other "structural" models in econometrics. Newey, Powell and Vella [13] investigated estimation and inference with a triangular-array version of (1.1). In that setup, equations relate $X_i$ and $W_i$, and the disturbances of these equations are connected to $U_i$. Newey and Powell [12] proposed a series estimator for $g$ in (1.1) and gave sufficient conditions for its consistency, but did not obtain a rate of convergence. Darolles, Florens and Renault [3] developed a kernel estimator for a special case of (1.1) and obtained its rate of convergence. This rate is slower than that obtained here. However, Darolles, Florens and Renault [3] make assumptions that conflict with ours, and it is not known whether their rate is optimal under their assumptions.

Further related work on inverse problems includes that of Wahba [17], Tikhonov and Arsenin [15], Groetsch [7], Nashed and Wahba [11] and Van Rooij and Ruymgaart [16].

We shall give a relatively detailed treatment, together with proofs, of results in the case where the instrumental variable is univariate. This setting is arguably of greatest interest to statisticians. Extensions to multivariate cases will be outlined.

## 2. Model and estimators in bivariate case.

2.1. *Model.* Let $(U_i, W_i, X_i, Y_i)$, for $i \geq 1$, be independent and identically distributed 4-vectors, and assume they follow a model satisfying (1.1) and (1.2). We shall suppose that $(W_i, X_i, Y_i)$, for $1 \leq i \leq n$, are observed, and that the distribution of $(X_i, W_i)$ is confined to the unit square.

Denote by $f_X$, $f_W$ and $f_{XW}$ the marginal densities of $X$ and $W$, and the joint density of $X$ and $W$, respectively, and define the linear operator $T$ on the space of square-integrable functions on $[0,1]^2$ by

$$(T\psi)(z) = \int t(x,z)\psi(x)\,dx,$$

where

$$t(x,z) = \int f_{XW}(x,w) f_{XW}(z,w)\,dw.$$

The following assumption characterizes the strength of association we require between $X$ and $W$:

(2.1) $\qquad\qquad\qquad T$ is nonsingular.

To appreciate the nature of (2.1), observe that if $X$ and $W$ are independent, then $T$ maps each function $\psi$ to a constant multiple of $f_X$, and so (2.1) fails. However, if (2.1) holds, then since it may be proved from (1.1) and (1.2) that

(2.2) $\qquad\qquad E_W\{E(Y|W) f_{XW}(z, W)\} = (Tg)(z),$

$g$ may be recovered by inversion of $T$,

(2.3) $\qquad\qquad g(z) = E_W\{E(Y|W)(T^{-1} f_{XW})(z, W)\}.$

This property suggests an estimator, which we shall develop in Section 2.2.

Observe that (2.2) is a Fredholm equation of the first kind, and generates an ill-posed inverse problem if, as is usually the case, zero is a limit point of the eigenvalues of $T$. In that case, $T^{-1}$ is not a bounded, continuous operator. For the purpose of estimation, we shall deal with this problem in Section 2.2 by replacing $T^{-1}$ by $(T + a_n)^{-1}$, where $a_n$ is a positive ridge parameter converging to zero as $n \to \infty$.

2.2. *Generalized kernel estimator.* Let $f_{XW}$ have $r$ continuous derivatives with respect to any combination of its arguments. Let $K_h(\cdot, \cdot)$ denote a generalized kernel function, with the properties $K_h(u,t) = 0$ if $u > t$ or $u < t - 1$,

for all $t \in [0, 1]$

(2.4)
$$h^{-(j+1)} \int_{t-1}^{t} u^j K_h(u,t)\, du = \begin{cases} 1, & \text{if } j = 0, \\ 0, & \text{if } 1 \leq j \leq r - 1. \end{cases}$$

Here, $h > 0$ denotes a bandwidth, and the kernel is considered in generalized form only to overcome edge effects. In particular, if $h$ is small and $t$ is not close to either 0 or 1, then we may take $K_h(u,t) = K(u/h)$, where $K$ is an $r$th order kernel. If $t$ is close to 1, then we may take $K_h(u,t) = L(u/h)$, where $L$ is a bounded, compactly supported function satisfying

$$\int_0^\infty u^j L(u)\, du = \begin{cases} 1, & \text{if } j = 0, \\ 0, & \text{if } 1 \leq j \leq r - 1. \end{cases}$$

And if $t$ is close to 0, then we may take $K_h(u,t) = L(-u/h)$. There are, of course, other ways of overcoming the edge-effect problem, but the "boundary kernel" approach above is also appropriate.

We require two estimators of $f_{XW}$, the second a leave-one-out estimator,

$$\hat{f}_{XW}(x,w) = \frac{1}{nh^2} \sum_{i=1}^{n} K_h(x - X_i, x) K_h(w - W_i, w),$$

$$\hat{f}_{XW}^{(-i)}(x,w) = \frac{1}{(n-1)h^2} \sum_{1 \leq j \leq n\,:\, j \neq i} K_h(x - X_j, x) K_h(w - W_j, w).$$

We use $\hat{f}_{XW}$ to construct the following estimators of $t(x,z)$ and the transformation $T$:

$$\hat{t}(x,z) = \int \hat{f}_{XW}(x,w) \hat{f}_{XW}(z,w)\, dw, \qquad (\widehat{T}\psi)(z) = \int \hat{t}(x,z) \psi(x)\, dx.$$

Let $a_n > 0$; we shall use it as a ridge parameter when inverting $\widehat{T}$, defining $\widehat{T}^+ = (\widehat{T} + a_n I)^{-1}$, where $I$ is the identity operator. Reflecting (2.3), our

estimator of $g$ is

$$\hat{g}(x) = \frac{1}{n}\sum_{i=1}^{n}(\widehat{T}^{+}\hat{f}_{XW}^{(-i)})(z, W_i)Y_i.$$

An alternative approach would be to develop a spectral expansion of $\widehat{T}$, truncate it to a finite series, and invert this series. The smoothing parameter now becomes the number of terms in the series, rather than the ridge, $a_n$. Theory may be developed for this "spectral cut-off" approach, too. However, it appears to require regularity conditions on spacings between adjacent eigenvalues of $T$, as well as a condition on their rate of decrease (see A.3 in Section 4.1), and for this reason we do not pursue it here.

2.3. *Orthogonal series estimator.* This technique is based on empirically transforming the marginal distributions of $W$ and $X$ to uniform, and exploiting the relatively simple character of the problem in that case. To appreciate this point, assume for the time being that both marginals are in fact uniform on $[0,1]$, and let $\chi_1, \chi_2, \ldots$ denote an orthonormal basis for $L_2[0,1]$. In practice, one would usually take $\{\chi_j\}$ to be the cosine sequence, although there are many other options.

Let $f_{XW}(x,w) = \sum_j \sum_k q_{jk}\chi_j(x)\chi_k(w)$ denote the generalized Fourier expansion of $f_{XW}$, and put $Q = (q_{jk})$, $p_j = E\{Y\chi_j(W)\}$, $\gamma_j = E\{g(X)\chi_j(X)\}$, $p = (p_j)$ and $\gamma = (\gamma_j)$, the latter two quantities being column vectors. By (1.1) and (1.2), $QQ'\gamma = Qp$ and, therefore, $\gamma = (QQ')^{-1}Qp$. [This is really another way of writing (2.3); observe that the operator $T$ takes $g$ to a function of which the $j$th Fourier coefficient is $(QQ'\gamma)_j$.] Hence, the problem of estimating the Fourier coefficients $\gamma_j$ of $g$ reduces to one of estimating $p_j$ and $q_{jk}$.

Next we describe how to solve the latter problem in general cases, where marginal distributions are not uniform. First transform the marginals, by computing $\widehat{W}_i = \widehat{F}_W(W_i)$ and $\widehat{X}_i = \widehat{F}_W(X_i)$, where $\widehat{F}_W$ and $\widehat{F}_X$ denote the empirical distribution functions of the data $W_1, \ldots, W_n$ and $X_1, \ldots, X_n$, respectively. Put $\hat{q}_{jk} = n^{-1}\sum_i \chi_j(\widehat{W}_i)\chi_k(\widehat{X}_i)$ and $\hat{p}_j = n^{-1}\sum_i \chi_j(\widehat{W}_i)Y_i$. Let $\widehat{Q}$ be the $m \times m$ matrix that has $\hat{q}_{jk}$ in position $(j,k)$, and set

$$\widehat{\gamma} = (\widehat{\gamma}_j) = (\widehat{Q}\widehat{Q}' + a_n I_m)^{-1}\widehat{Q}\hat{p},$$

where $a_n$ denotes a ridge parameter and $I_m$ is the $m \times m$ identity. Our estimator of $g$ is

$$\bar{g}(x) = \sum_{j=1}^{m}\widehat{\gamma}_j\chi_j(x).$$

In this estimator the number of terms, $m$, in the approximating Fourier series is the main smoothing parameter. It is relatively awkward to derive theory for the orthogonal series method, owing to the fact that the transformed data $\widehat{W}_i$ and $\widehat{X}_i$ are not independent, and to the difficulty of dealing theoretically with the large random matrix $\widehat{Q}$. Nevertheless, we shall show in Section 4 that, under restrictions, the orthogonal series technique has optimal performance.

**3. Model and estimators in the multivariate case.** In the model at (1.1) the explanatory variable $X$ is endogenous, that is, determined within the model. When the model is multivariate, there is an opportunity for dividing the explanatory variable, which is now a vector, into two parts, one endogenous and the other determined outside the model, or exogenous.

We take $(Y, X, Z, W, U)$ to be a vector, where $Y$ and $U$ are scalars, $X$ and $W$ are supported on $[0,1]^p$, and $Z$ is supported on $[0,1]^q$. Generalizing (1.1) and (1.2), the model is

$$Y_i = g(X_i, Z_i) + U_i, \qquad E(U_i | Z_i, W_i) = 0,$$

where $(Y_i, X_i, Z_i, W_i, U_i)$, for $i \geq 1$, are independent and identically distributed as $(Y, X, Z, W, U)$. Thus, $X$ and $Z$ are endogenous and exogenous explanatory variables, respectively. Data $(Y_i, X_i, Z_i, W_i)$, for $1 \leq i \leq n$, are observed.

Let $f_{XZW}$ denote the density of $(X, Z, W)$, write $f_Z$ for the density of $Z$, and for each $x_1, x_2 \in [0,1]^p$ put

$$t_z(x_1, x_2) = \int f_{XZW}(x_1, z, w) f_{XZW}(x_2, z, w)\,dw,$$

the analogue of $t(x_1, x_2)$ in Section 2. Define the operator $T_z$ on $L_2[0,1]^p$ by

$$(T_z \psi)(x) = \int t_z(\xi, x) \psi(\xi)\,d\xi.$$

Analogously to (2.3), it may be proved that, for each $z$ for which $T_z^{-1}$ exists,

$$g(x, z) = f_Z(z) E_{W|Z}\{E(Y | Z = z, W)(T_z^{-1} f_{XZW})(x, z, W) | Z = z\},$$

where $E_{W|Z}$ denotes the expectation operator with respect to the distribution of $W$ conditional on $Z$. In this formulation, $(T_z^{-1} f_{XZW})(x, z, W)$ denotes the result of applying $T_z^{-1}$ to the function $f_{XZW}(\cdot, z, W)$ and evaluating the resulting function at $x$.

To construct an estimator of $g(x, z)$, given $h > 0$ and $p$-vectors $x = (x^{(1)}, \ldots, x^{(p)})$ and $\xi = (\xi^{(1)}, \ldots, \xi^{(p)})$, let $K_{p,h}(x, \xi) = \prod_{1 \leq j \leq p} K_h(x^{(j)}, \xi^{(j)})$, put $K_{q,h}(z, \zeta)$ analogously for $q$-vectors $z$ and $\zeta$, let $h_x, h_z > 0$, and define

$$\hat{f}_{XZW}(x, z, w) = \frac{1}{n h_x^{2p} h_z^q} \sum_{i=1}^{n} K_{p,h_x}(x - X_i, x)$$

$$\times K_{q,h_z}(z - Z_i, z) K_{p,h_x}(w - W_i, w),$$

$$\hat{f}_{XZW}^{(-i)}(x,z,w) = \frac{1}{(n-1)h_x^{2p}h_z^q} \sum_{1 \leq j \leq n: j \neq i} K_{p,h_x}(x - X_j, x)$$

$$\times K_{q,h_z}(z - Z_j, z) K_{p,h_x}(w - W_j, w),$$

$$\hat{t}_z(x_1, x_2) = \int \hat{f}_{XZW}(x_1, z, w) \hat{f}_{XZW}(x_2, z, w) \, dw$$

and

$$(\widehat{T}_z \psi)(x, z, w) = \int \hat{t}_z(\xi, x) \psi(\xi, z, w) \, d\xi,$$

where $\psi$ is a function from $\mathbb{R}^{2p+q}$ to the real line. Then the estimator of $g(x, z)$ is

$$\hat{g}(x, z) = \frac{1}{n} \sum_{i=1}^{n} (\widehat{T}_z^+ \hat{f}_{XZW}^{(-i)})(x, z, W_i) Y_i K_{q, h_z}(z - Z_i, z).$$

## 4. Theoretical properties.

4.1. *Kernel method for bivariate case.* The invertibility of $T$ is central to our ability to successfully resolve $g$ from data, and so it comes as no surprise to find that rates of convergence of estimators of $g$ hinge on the rate at which the eigenvalues of $T$, say $\lambda_1 \geq \lambda_2 \geq \cdots > 0$, converge to 0. Therefore, our regularity conditions will be framed in terms of an eigenexpansion representation of $T$. To this end, let $\phi_j$ denote an eigenfunction of $T$ with eigenvalue $\lambda_j$, normalized so that $\phi_1, \phi_2, \ldots$ is an orthonormal basis for the space of square-integrable functions on the interval $[0, 1]$. Then we may write

$$t(x, z) = \sum_{j=1}^{\infty} \lambda_j \phi_j(x) \phi_j(z),$$

(4.1)
$$f_{XW}(x, z) = \sum_{j=1}^{\infty} \sum_{k=1}^{\infty} d_{jk} \phi_j(x) \phi_k(z),$$

$$g(x) = \sum_{j=1}^{\infty} b_j \phi_j(x),$$

where $d_{jk}$ and $b_j$ denote generalized Fourier coefficients of $f_{XW}$ and $g$, respectively.

Next we state regularity conditions. Assumption A.1 is equivalent to the intersection of (1.1) and (1.2); A.3 gives smoothness conditions, expressed through the eigen-expansion of $T$; A.2 and A.3 together imply that $T$ is a bounded Hilbert–Schmidt operator and, hence, compact; and A.4 describes the sizes of tuning parameters. The invertibility condition (2.1) is equivalent to asking that each $\lambda_j > 0$, which in turn implies part of A.3.

Below, in condition A.3, we shall introduce constants $\alpha, \beta > 0$, for which

$$A_1 \equiv \max\left(\frac{2\alpha + 2\beta - 1}{2\beta - \alpha}, \frac{5}{2}\frac{2\alpha + 2\beta - 1}{4\beta - \alpha + 1}, 2\right) > 0,$$

$$0 < A_2 \equiv \frac{1}{2r}\frac{2\alpha + 2\beta - 1}{2\beta + \alpha} \leq A_3 \equiv \min\left\{\frac{1}{2}\frac{2\beta - \alpha}{2\beta + \alpha}, \frac{4\beta - \alpha + 1}{5(2\beta + \alpha)}\right\}.$$

Therefore, it is possible to choose an integer $r \geq A_1$ and a constant $\gamma \in [A_2, A_3]$; such values will be used below. Let $C > 0$ be an arbitrarily large but fixed constant, let $\alpha, \beta > 0$, and denote by $\mathcal{G} = \mathcal{G}(C, \alpha, \beta)$ the class of distributions $G$ of $(X, W, Y)$ that satisfy A.1–A.3 below.

Regarding the smoothness assumed of $f_{XW}$ in A.2, we mention that our minimax rates do not alter if $f_{XW}$ is smoother than specified. The rates are optimized for smoothness of $g$, given enough smoothness of $f_{XW}$. In condition A.3, the lower bound on $\alpha$ seems difficult to relax and, in fact, it has close analogues in related contexts, for example, in work on convergence rates in functional data problems.

The upper bound on $\alpha$, however, seems more likely to be tied to our method of proof. One approach to relaxing the bound might be to draw inspiration from a modified approach to Tikhonov regularization (see [10]) and use, as the ridged inverse, $(T + a_n D^{2\beta - 1})^{-1}$ rather than $(T + a_n I)^{-1}$. Here, if $2\beta - 1$ were an integer, $D^{2\beta - 1}$ would denote the $(\beta - 1)$st power of the differential operator; if $2\beta - 1$ were strictly greater than its integer part, $\ell$ say, then $D^{2\beta - 1}$ would involve taking the convolution of $g^{(\ell)}(t) - g^{(\ell)}(0)$ against the kernel $|t|^{\ell - 2\beta}$. However, this approach requires a direct relationship between the smoothness of $g$, as expressed through the size of $\beta$ in the formula $|b_j| \leq Cj^{-\beta}$, and its smoothness in the more conventional sense of differentiation. We have avoided making assumptions about this relationship. In particular, as our results are presently formulated, $g$ does not need to be continuous, let alone differentiable, no matter how large or small $\beta$ might be.

A.1. The data $(X_i, W_i, Y_i)$ are independent and identically distributed as $(X, W, Y)$, where $(X, W)$ is supported on $[0, 1]^2$ and $E\{Y - g(X)|W = w\} \equiv 0$.

A.2. The distribution of $(X, W)$ has a density, $f_{XW}$, with $r$ derivatives (when viewed as a function restricted to $[0, 1]^2$) bounded uniformly in absolute value by $C$; and the functions $E(Y^2|W = w)$ and $E(Y^2|X = x, W = w)$ are bounded uniformly by $C$.

A.3. The constants $\alpha$ and $\beta$ satisfy $\alpha > 1$, $\beta > \frac{1}{2}$ and $\beta - \frac{1}{2} \leq \alpha < 2\beta$. Moreover, $|b_j| \leq Cj^{-\beta}$, $j^{-\alpha} \leq C\lambda_j$ and $\sum_{k\geq 1} |d_{jk}| \leq Cj^{-\alpha/2}$ for all $j \geq 1$.

A.4. The parameters $a_n$ and $h$ satisfy $a_n \asymp n^{-\alpha/(2\beta+\alpha)}$ and $h \asymp n^{-\gamma}$ as $n \to \infty$, where $c_n \asymp d_n$ for positive constants $c_n$ and $d_n$ means that $c_n/d_n$ is bounded away from zero and infinity.

A.5. The function $K_h(\cdot, \cdot)$ satisfies (2.4); for each $t \in [0,1]$, $K_h(h\cdot, t)$ is supported on $[(t-1)/h, t/h] \cap \mathcal{K}$, where $\mathcal{K}$ is a compact interval not depending on $t$; and

$$\sup_{h>0, t\in[0,1], u\in\mathcal{K}} |K_h(hu,t)| < \infty.$$

THEOREM 4.1. *As $n \to \infty$,*

$$\sup_{G\in\mathcal{G}} \int_0^1 E_G\{\hat{g}(t) - g(t)\}^2 \, dt = O(n^{-(2\beta-1)/(2\beta+\alpha)}).$$

More generally, it may be proved that if a particular distribution of $(X, W, Y)$ satisfies A.1, and if $E(Y^2) < \infty$ and the density $f_{XW}$ is continuous on $[0,1]$, then $a_n$ and $h$ can be chosen so that $\int E_G(\hat{g} - g)^2 \to 0$ as $n \to \infty$. Similar results, guaranteeing consistent estimation but without a convergence rate, may be derived in the settings of Sections 4.2 and 4.3.

4.2. *Orthogonal series method for bivariate case.* We shall simplify theory by assuming the Fourier coefficients $q_{jk}$ satisfy a strong diagonality condition. Under this assumption it is sufficient to work with a strongly diagonal form of $\widehat{Q}$, where we redefine $\hat{q}_{jk} = 0$ if $|j - k| \geq N$ (where $N$ is permitted to increase slowly with $n$), and leave $\hat{q}_{jk}$ unchanged otherwise. With this alteration to $\hat{q}_{jk}$, let $\widehat{Q} = (\hat{q}_{jk})$ be the indicated $m \times m$ matrix.

Recall from Section 2.3 that $\chi_1, \chi_2, \ldots$ is an orthonormal basis for $L_2[0,1]$. Let $F_W$ and $F_X$ denote the marginal distribution functions of $W$ and $X$, put $\widetilde{W} = F_W(W)$ and $\widetilde{X} = F_X(X)$, and let $f_{\widetilde{W}\widetilde{X}}$ denote the joint density of $(\widetilde{W}, \widetilde{X})$. Write $f_{\widetilde{W}\widetilde{X}}(w,x) = \sum_j \sum_k q_{jk}\chi_j(x)\chi_k(w)$ and $g(x) = \sum_j \gamma_j \chi_j(x)$ for the generalized Fourier transforms of these functions. Recall that we require the transformation represented by $QQ'$ to be invertible, so we may define $Q^{-1} = (q_{jk}^{(-1)})$ to be a generalized inverse of $Q$.

Given constants $\alpha \geq 2$, $\beta \geq \frac{1}{2}$ and $C_1, C_2 > 0$, let $\mathcal{H} = \mathcal{H}(C_1, C_2, \alpha, \beta)$ denote the class of distributions $G$ of $(\widetilde{W}, \widetilde{X}, Y)$ for which

$$E\{Y - g(\widetilde{X})|\widetilde{W} = w\} \equiv 0, \qquad |q_{jk}| \leq C_1\{\max(j,k)\}^{-\alpha/2}\exp(-C_2|j-k|),$$

$$|q_{jk}^{(-1)}| \leq C_1\{\max(j,k)\}^{\alpha/2}\exp(-C_2|j-k|),$$

$$|p_j| \leq C_1 j^{-\beta}, \qquad E(Y^4) < C_1,$$

where the bounds are assumed to hold uniformly in $1 \leq j, k < \infty$.

THEOREM 4.2. *Let $\{\chi_j\}$ denote the orthonormalized version of the cosine series on $[0,1]$. Take $\alpha \geq 2$ and $\beta \geq \frac{1}{2}$, and assume $a_n \asymp m^{-\alpha}$, $m \asymp n^{1/(2\beta+\alpha)}$, $N/\log n \to \infty$ and $N = O(n^\varepsilon)$ for all $\varepsilon > 0$. Then, as $n \to \infty$,*

$$\sup_{G \in \mathcal{H}} \int_0^1 E_G(\bar{g} - g)^2 = O(n^{-(2\beta-1)/(2\beta+\alpha)}).$$

4.3. *Kernel method for multivariate case.* For each $z \in [0,1]^q$, let $\{\phi_{z1}, \phi_{z2}, \ldots\}$ denote the orthonormalized sequence of eigenvectors, and $\lambda_{z1} \geq \lambda_{z2} \geq \cdots > 0$ the respective eigenvalues of the operator $T_z$. Assume that $\{\phi_{zj}\}$ forms an orthonormal basis of $L_2[0,1]^p$. Analogously to (4.1),

$$t_z(x_1, x_2) = \sum_{j=1}^\infty \lambda_{zj} \phi_{zj}(x_1) \phi_{zj}(x_2),$$

$$f_{XZW}(x, z, w) = \sum_{j=1}^\infty \sum_{k=1}^\infty d_{zjk} \phi_{zj}(x) \phi_{zk}(z),$$

$$g(x, z) = \sum_{j=1}^\infty b_{zj} \phi_{zj}(x),$$

where the $d_{zjk}$'s and $b_{zj}$'s are generalized Fourier coefficients.

Put $\tau = 2r/(2r+q)$. If $\alpha, \beta > 0$ denote constants satisfying MV.3 below, then

$$B_1 \equiv \max\left\{p\frac{2\alpha + 2\beta - 1}{2\beta - \alpha},\right.$$

$$\left.\left(\frac{2}{5p}\frac{2\alpha + 10\beta + 1}{2\beta + \alpha} - \frac{6}{5p}\right)^{-1}\left(\frac{2\alpha + 2\beta - 1}{2\beta + \alpha} + \frac{3q}{5p}\right), 2\right\} > 0,$$

$$0 < B_2 \equiv \frac{\tau}{2r}\frac{2\alpha + 2\beta - 1}{2\beta + \alpha} \leq B_3 \equiv \min\left\{\frac{\tau}{2p}\frac{2\beta - \alpha}{2\beta + \alpha}, \frac{1}{5p}\left(\tau\frac{10\beta + 2\alpha}{2\beta + \alpha} - 3\right)\right\}.$$

Choose $r \geq B_1$ and $\gamma \in [B_2, B_3]$. We make the following assumptions, of which the first five are respectively analogous to A.1–A.5 in Section 4.1. Let $C > 0$.

MV.1. The data $(X_i, W_i, Z_i, Y_i)$ are independent and identically distributed as $(X, W, Z, Y)$, where $X$, $W$ and $Z$ are supported on $[0,1]^p$, $[0,1]^p$ and $[0,1]^q$, respectively, and $E\{Y - g(X, Z)|Z = z, W = w\} \equiv 0$.

MV.2. The distribution of $(X, Z, W)$ has a density, $f_{XZW}$, with $r$ derivatives of all types (when viewed as a function restricted to $[0,1]^{2p+q}$), each derivative bounded in absolute value by $C$; $g(x, z)$ and $b_{zj}$ have $r$ partial derivatives with respect to $z$, bounded in absolute value by $C$, uniformly in $x$ and $z$; and the functions $E(Y^2|Z = z, W = w)$ and $E(Y^2|X = x, Z = z, W = w)$ are bounded uniformly by $C$.

MV.3. The constants $\alpha, \beta$ satisfy $\alpha > 1$, $\beta > \frac{1}{2}$ and $\beta - \frac{1}{2} \leq \alpha < 2\beta$. Moreover, $|b_{zj}| \leq Cj^{-\beta}$, $j^{-\alpha} \leq C\lambda_{zj}$ and $\sum_{k\geq 1} |d_{zjk}| \leq Cj^{-\alpha/2}$, uniformly in $z \in [0,1]^q$, for all $j \geq 1$.

MV.4. The parameters $a_n$, $h_x$ and $h_z$ satisfy $a_n \asymp n^{-\alpha\tau/(2\beta+\alpha)}$, $h \asymp n^{-\gamma}$, $h_z \asymp n^{-1/(2r+q)}$ as $n \to \infty$.

MV.5. The function $K_h(\cdot, \cdot)$ satisfies A.5.

MV.6. For each $z \in [0,1]^q$, the functions $\phi_{zj}$ form an orthonormal basis for $L_2[0,1]^p$, and $\sup_x \sup_z \max_j |\phi_{zj}(x)| < \infty$.

Let $\mathcal{M} = \mathcal{M}(C, \alpha, \beta)$ denote the class of distributions of $(X, W, Z, Y)$ that satisfy MV.1–MV.3 and MV.6.

THEOREM 4.3. *As $n \to \infty$,*

$$\sup_{G \in \mathcal{M}} \sup_{z \in [0,1]^q} \int_{[0,1]^p} E_G\{\hat{g}(x,z) - g(x,z)\}^2 \, dx = O(n^{-\tau(2\beta-1)/(2\beta+\alpha)}).$$

4.4. *Optimality.* The convergence rates expressed by Theorems 4.1–4.3 are optimal in those contexts, in a minimax sense. Indeed, let $\tilde{g}$ denote any measurable functional of that data which is itself a measurable function on $[0,1]$ (in the cases of Theorems 4.1 and 4.2) or on $[0,1]^p$ (in the setting of Theorem 4.3); let $\mathcal{C}$ denote $\mathcal{G}$, $\mathcal{H}$ or $\mathcal{M}$ in the cases of Theorems 4.1–4.3, respectively; and put $\tau = 1$ in the contexts of Theorems 4.1 and 4.2, and $\tau = 2r/(2r+q)$ for Theorem 4.3.

THEOREM 4.4.

(4.2) $$\liminf_{n\to\infty} n^{\tau(2\beta-1)/(2\beta+\alpha)} \inf_{\tilde{g}} \sup_{G \in \mathcal{C}} \int E_G(\tilde{g} - g)^2 > 0.$$

In the multivariate setting of Section 4.3 we interpret the integral at (4.2) as

$$\int_{[0,1]^p} E_G\{\tilde{g}(x,z) - g(x,z)\}^2 \, dx,$$

and interpret Theorem 4.4 as stating that, for this representation, (4.2) holds for each $z \in [0,1]^q$.

**5. Monte Carlo experiments.** This section reports the results of a Monte Carlo investigation of the finite-sample performance of the kernel estimator for the bivariate model. The estimator is the one described in Section 2.3, although our method is not optimized for theoretical performance. In particular, we took $K$ to be a second-order kernel.

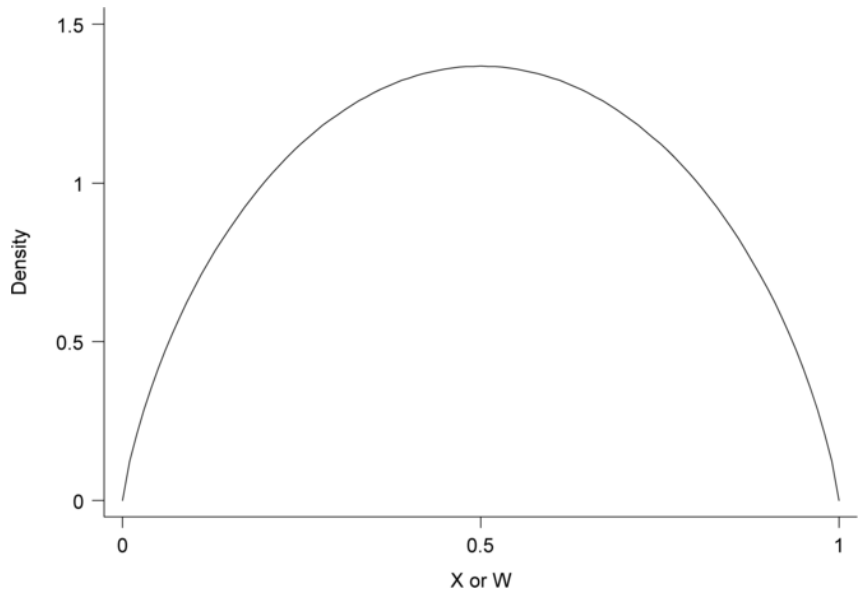

Fig. 1. *Density of X and W used in Monte Carlo experiments.*

Samples of size $n = 200$ were generated from the model determined by

$$f_{XW}(x,w) = 2C_f \sum_{j=1}^{\infty} (-1)^{j+1} j^{-1} \sin(j\pi x) \sin(j\pi w), \qquad 0 \leq x, w \leq 1;$$

$$g(x) = 2^{1/2} \sum_{j=1}^{\infty} (-1)^{j+1} j^{-2} \sin(j\pi x), \qquad Y = E\{g(X)|W = w\} + V,$$

where $C_f$ is a normalization constant and $V$ is distributed as Normal $N(0, 0.01)$. For computational purposes, the infinite series were truncated at $j = 100$. Figure 1 shows a graph of the marginal distributions of $X$ and $W$, which are identical. The solid line in Figure 2 depicts $g(x)$. The kernel function is the Epanechnikov kernel, $K(x) = 0.75(1 - x^2)$ for $|x| \leq 1$.

Each experiment consisted of estimating $g$ at the 19 points, $x = 0.05, 0.10, \ldots, 0.95$. The experiments were carried out in GAUSS using GAUSS pseudo-random number generators. There were 1000 Monte Carlo replications in each experiment.

Table 1 shows the performance of the estimator, $\hat{g}$, as a function of the bandwidth, $h$, and the ridge parameter, $a_n$. The quantities Bias$^2$, Var and MSE in the table were calculated as the averages, over the 19 values of $x$, of Monte Carlo approximations to pointwise squared bias, variance and mean squared error, respectively, at those points; the pointwise values were computed by averaging over the 1000 Monte Carlo simulations.

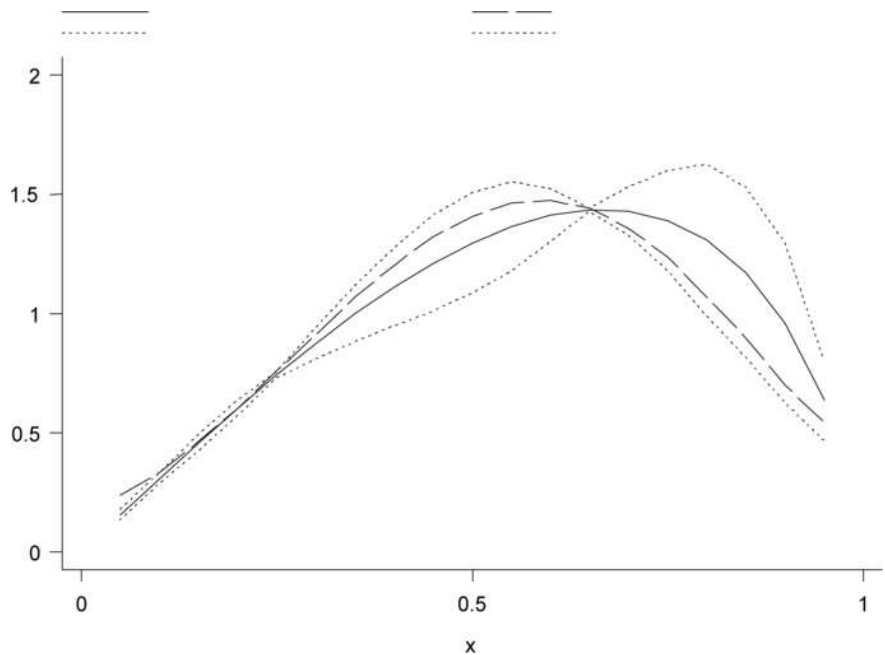

FIG. 2. *Graph of 95% estimation band. The solid, dashed and dotted lines show $g$, $E(\hat{g})$ and the 95% estimation band, respectively.*

TABLE 1
*Results of Monte Carlo experiments*

| $a_n$ | $h$ | Bias$^2$ | Var | MSE |
|---|---|---|---|---|
| 0.05 | 0.10 | 0.0039 | 0.0321 | 0.0361 |
|  | 0.20 | 0.0065 | 0.0162 | 0.0227 |
|  | 0.30 | 0.0262 | 0.0119 | 0.0381 |
|  | 0.40 | 0.0525 | 0.0087 | 0.0612 |
| 0.10 | 0.10 | 0.0118 | 0.0221 | 0.0339 |
|  | 0.20 | 0.0105 | 0.0115 | 0.0215 |
|  | 0.30 | 0.0141 | 0.0078 | 0.0219 |
|  | 0.40 | 0.0263 | 0.0062 | 0.0325 |
| 0.15 | 0.10 | 0.0224 | 0.0190 | 0.0414 |
|  | 0.20 | 0.0165 | 0.0098 | 0.0263 |
|  | 0.30 | 0.0149 | 0.0063 | 0.0212 |
|  | 0.40 | 0.0220 | 0.0049 | 0.0269 |
| 0.20 | 0.10 | 0.0335 | 0.0174 | 0.0508 |
|  | 0.20 | 0.0268 | 0.0081 | 0.0349 |
|  | 0.30 | 0.0214 | 0.0058 | 0.0272 |
|  | 0.40 | 0.0252 | 0.0044 | 0.0295 |

Results are illustrated graphically in Figure 2 for the case $h = 0.2$ and $a_n = 0.1$. The figure shows $g(x)$ (solid line), the Monte Carlo approximation to $E\{\hat{g}(x)\}$ (dashed line) and a 95% pointwise "estimation band." The band connects the points $g(x_j) \pm \delta_j$, for $j = 1, \ldots, 19$, where each $\delta_j$ is chosen so that the interval $[g(x_j) - \delta_j, g(x_j) + \delta_j]$ contains 95% of the 1000 simulated values of $\hat{g}(x_j)$. The figure shows, not surprisingly, that $\hat{g}$ is somewhat biased, but that the shape of $E\hat{g}$ is similar to that of $g$.

## 6. Technical arguments.

6.1. *Proof of Theorem* 4.1. (The "big oh" bounds that we shall derive below apply uniformly in $G \in \mathcal{G}$, although for the sake of simplicity we shall not make this qualification.) Put $T^+ = (T + a_n I)^{-1}$, let $\|\cdot\|$ denote the usual $L_2$ norm for functions from the interval $[0, 1]$ to the real line, and given a functional $\chi$ from $L_2[0, 1]$ to itself, set

$$\|\chi\| = \sup_{\psi \in L_2[0,1]\,:\,\|\psi\|=1} \|\chi(\psi)\|.$$

For future reference, we note that A.3 and A.4 imply that

(6.1) $\quad n^{\{1/(2\beta+\alpha)\}-1} a_n^{-1} + h^{2r} a_n^{-2} = O(n^{-(2\beta-1)/(2\beta+\alpha)}).$

Define

$$D_n(z) = \int g(x) f_{XW}(x, w) T^+ (\hat{f}_{XW} - f_{XW})(z, w)\, dx\, dw,$$

$$A_{n1}(z) = \frac{1}{n} \sum_{i=1}^n (T^+ f_{XW})(z, W_i) Y_i,$$

$$A_{n2}(z) = \frac{1}{n} \sum_{i=1}^n \{T^+ (\hat{f}_{XW}^{(-i)} - f_{XW})\}(z, W_i) Y_i - D_n(z),$$

$$A_{n3}(z) = \frac{1}{n} \sum_{i=1}^n \{(\widehat{T}^+ - T^+) f_{XW}\}(z, W_i) Y_i + D_n(z),$$

$$A_{n4}(z) = \frac{1}{n} \sum_{i=1}^n \{(\widehat{T}^+ - T^+)(\hat{f}_{XW}^{(-i)} - f_{XW})\}(z, W_i) Y_i.$$

Then $\hat{g} = A_{n1} + \cdots + A_{n4}$, and so the theorem will follow if we prove that

(6.2) $\quad E\|A_{n1} - g\|^2 = O(n^{-(2\beta-1)/(2\beta+\alpha)}),$

(6.3) $\quad E\|A_{nj}\|^2 = O(n^{-(2\beta-1)/(2\beta+\alpha)}) \quad$ for $j = 2, 3, 4$.

To derive (6.2), note that $EA_{n1} - g = -a_n \sum_{j \geq 1} b_j (\lambda_j + a_n)^{-1} \phi_j$. Therefore,

$$\|EA_{n1} - g\|^2 = a_n^2 \sum_{j=1}^\infty \frac{b_j^2}{(\lambda_j + a_n)^2}.$$

Divide the last-written series up into the sum over $j \leq J \equiv a_n^{-1/\alpha}$, and the complementary part, thereby bounding the right-hand side by $a_n^2 \sum_{j \leq J} (b_j/\lambda_j)^2 + \sum_{j > J} b_j^2$; and use A.3 and A.4 to bound each of these terms, hence, proving that

(6.4) $$\|EA_{n1} - g\|^2 = O(n^{-(2\beta-1)/(2\beta+\alpha)}).$$

Using A.2, we deduce that

$$n \operatorname{var}\{A_{n1}(z)\} \leq E[\{(T^+ f_{XW})(z,W)Y\}^2]$$
$$= E[\{(T^+ f_{XW})(z,W)\}^2 E(Y^2|W)]$$
$$\leq \operatorname{const.} B_n,$$

where $B_n = E[\{(T^+ f_{XW})(z,W)\}^2]$ and, here and below, "const." will denote a positive constant, different at different appearances. It can be proved, from an expansion of $T^+ f_{XW}(z,w)$ in its generalized Fourier series, that

$$B_n = \sum_{j=1}^{\infty} \sum_{k=1}^{\infty} \sum_{\ell=1}^{\infty} \frac{d_{jk} d_{j\ell}}{(\lambda_j + a_n)^2} E\{\phi_k(W)\phi_\ell(W)\}$$
$$\leq \operatorname{const.} \sum_{j=1}^{\infty} \sum_{k=1}^{\infty} \sum_{\ell=1}^{\infty} \frac{|d_{jk} d_{j\ell}|}{(\lambda_j + a_n)^2}$$
$$\leq \operatorname{const.} \sum_{j=1}^{\infty} \frac{\lambda_j}{(\lambda_j + a_n)^2}.$$

Therefore,

$$\int \operatorname{var}\{A_{n1}(z)\} dz \leq \operatorname{const.} \frac{1}{n} \sum_{j=1}^{\infty} \frac{\lambda_j}{(\lambda_j + a_n)^2}.$$

From this point, using the argument leading to (6.4), we may prove that

$$E\|A_{n1} - EA_{n1}\|^2 = \int \operatorname{var}\{A_{n1}(z)\} dz$$
$$= O(n^{-1} a_n^{-(\alpha+1)/\alpha})$$
$$= O(n^{-(2\beta-1)/(2\beta+\alpha)}).$$

Result (6.2) is implied by this bound and (6.4).

Next we derive (6.3) in the case $j = 2$. Here and below, given a bivariate function $\phi(z,w)$, put $\phi_w(z) = \phi(z,w)$ and define $T^+\phi(z,w) = (T^+\phi_w)(z)$. Let

$$D_{ni}(z) = \int g(x) f_{XW}(x,w) T^+(\hat{f}_{XW}^{(-i)} - f_{XW})(z,w) \, dx \, dw,$$

$$A_{n21}(z) = \frac{1}{n}\sum_{i=1}^{n}\{T^+(\hat{f}_{XW}^{(-i)} - f_{XW})(z, W_i)Y_i - D_{ni}(z)\},$$

$$A_{n22}(z) = \frac{1}{n}\sum_{i=1}^{n}\{D_{ni}(z) - D_n(z)\},$$

in which notation $A_{n2} = A_{n21} + A_{n22}$. Write $\int A_{n21}(z)^2\,dz$ as a double series, and take the expected values of the terms one by one. It may be shown by tedious calculation that the total contribution of the terms equals $O\{h^{2r}(na_n^2)^{-1} + (nha_n)^{-2}\}$. Therefore,

(6.5) $\quad E\|A_{n21}\|^2 = O\{h^{2r}(na_n^2)^{-1} + (nha_n)^{-2}\} = O(n^{-(2\beta-1)/(2\beta+\alpha)}),$

where we used (6.1) to obtain the second identity. Furthermore,

$$A_{n22}(z) = -n^{-1}\int g(x)f_{XW}(x,w)T^+\hat{f}_{XW}(z,w)\,dx\,dw,$$

from which, noting (6.1), it may be deduced that

$$E\|A_{n22}\|^2 \leq \text{const.}(na_n)^{-2}E\left(\int |gf_{XW}\hat{f}|\right)^2$$
$$= O\{(na_n)^{-2}\} = O(n^{-(2\beta-1)/(2\beta+\alpha)}).$$

Property (6.3), in the case $j=2$, follows from this result and (6.5).

Next we derive (6.3) for $j=3$. Define $\Delta = \hat{T} - T$, an operator, and put

$$A_{n31} = -(I+T^+\Delta)^{-1}T^+\Delta g + D_n, \qquad A_{n32} = -(I+T^+\Delta)^{-1}T^+\Delta(A_{n1} - g).$$

Noting that $\hat{T}^+ - T^+ = -(I+T^+\Delta)^{-1}T^+\Delta T^+$, it can be seen that $A_{n3} = A_{n31} + A_{n32}$.

Let $\delta = h^{2r} + (nh)^{-1}$. Using standard, but tedious, moment calculations, it may be proved that $E(\hat{t}-t)^{2k} = O(\delta^k)$ for each integer $k \geq 1$, uniformly in the argument of $\hat{t}-t$. [The quantity $\delta$ involves $(nh)^{-1}$, rather than $(nh^2)^{-1}$, since the integral in the definition of $\hat{t}$ effectively removes one of the factors $h^{-1}$.] Therefore, since $\|\Delta\|^2 = \int(\hat{t}-t)^2$, then for each integer $k \geq 1$,

(6.6) $\qquad\qquad E\|\Delta\|^{2k} = O(\delta^k).$

At the end of this proof we shall show that, for each $k \geq 1$,

(6.7) $\qquad\qquad E\{\|(I+T^+\Delta)^{-1}\|^k\} = O(1)$

as $n \to \infty$. Hence, using the Cauchy–Schwarz inequality,

(6.8) $\quad \{E\|(I+T^+\Delta)^{-1}T^+\Delta\|^4\}^2 \leq E\|(I+T^+\Delta)^{-1}\|^8\|T^+\|^8 E\|\Delta\|^8$
$$= O(\delta^4/a_n^8).$$

From this result, and the Cauchy–Schwarz inequality again, we obtain
$$E\|A_{n32}\|^2 \leq \{E\|(I+T^+\Delta)^{-1}T^+\Delta\|^4 E\|A_{n1}-g\|^4\}^{1/2}$$
(6.9)
$$= O\{(\delta/a_n^2)^2 (E\|A_{n1}-g\|^4)^{1/2}\}$$
$$= O(n^{-(2\beta-1)/(2\beta+\alpha)}),$$
the final identity following using an argument similar to that leading to (6.2).

Put
$$B_{n1}(z) = \int \{\hat{f}_{XW}(x,w) - f_{XW}(x,w)\} f_{XW}(z,w) g(x) \, dx \, dw,$$

$$B_{n2}(z) = \int \{\hat{f}_{XW}(z,w) - f_{XW}(z,w)\} f_{XW}(x,w) g(x) \, dx \, dw,$$

$$B_{n3}(z) = \int \{\hat{f}_{XW}(x,w) - f_{XW}(x,w)\}\{\hat{f}_{XW}(z,w) - f_{XW}(x,w)\} g(x) \, dx \, dw,$$

$$B_{n11}(z) = \int \{E\hat{f}_{XW}(x,w) - f_{XW}(x,w)\} f_{XW}(z,w) g(x) \, dx \, dw,$$

$$B_{n12}(z) = \int \{\hat{f}_{XW}(x,w) - E\hat{f}_{XW}(x,w)\} f_{XW}(z,w) g(x) \, dx \, dw,$$

$$B_{n21}(z) = \int \{E\hat{f}_{XW}(z,w) - f_{XW}(z,w)\} f_{XW}(x,w) g(x) \, dx \, dw,$$

$$B_{n22}(z) = \int \{\hat{f}_{XW}(z,w) - E\hat{f}_{XW}(z,w)\} f_{XW}(x,w) g(x) \, dx \, dw.$$

In this notation, $\Delta g = B_{n1} + B_{n2} + B_{n3}$, $B_{n1} = B_{n11} + B_{n12}$, $B_{n2} = B_{n21} + B_{n22}$ and $T^+ B_{n2} = D_n$, whence
$$A_{n31} = -(I+T^+\Delta)^{-1}T^+(B_{n11}+B_{n12}+B_{n3})$$
$$+ (I+T^+\Delta)^{-1}T^+\Delta T^+(B_{n21}+B_{n22}).$$

Define
$$\tilde{A}_{n31} = -(I+T^+\Delta)^{-1}T^+(B_{n11}+B_{n12}+B_{n3}) + (I+T^+\Delta)^{-1}T^+\Delta T^+ B_{n21}.$$
Then
$$E\|A_{n31}\|^2 \leq \text{const.}\{E\|\tilde{A}_{n31}\|^2 + E\|(I+T^+\Delta)^{-1}T^+\Delta T^+ B_{n22}\|^2\}.$$
By (6.7) and the Cauchy–Schwarz inequality,
(6.10)
$$E\|\tilde{A}_{n31}\|^2 \leq \text{const.}(\|T^+ B_{n11}\|^4 + E\|T^+ B_{n12}\|^4$$
$$+ E\|T^+\Delta T^+ B_{n21}\|^4 + E\|T^+ B_{n3}\|^4)^{1/2}.$$

Since $\|B_{n11}\| + \|B_{n21}\| = O(h^r)$ and $\|T^+\| = O(a_n^{-1})$, then, by (6.1),
(6.11)
$$\|T^+ B_{n11}\| + \|T^+ B_{n21}\| \leq \|T^+\|(\|B_{n11}\| + \|B_{n21}\|)$$
$$= O(h^r a_n^{-1}) = O(n^{-(2\beta-1)/\{2(2\beta+\alpha)\}}).$$

Furthermore, with

$$\Delta_{jk} = \int \{\hat{f}_{XW}(x,w) - E\hat{f}_{XW}(x,w)\}\phi_j(x)\phi_k(x)\,dx\,dw,$$

we have

$$T^+ B_{n12}(z) = \sum_{j=1}^{\infty}\sum_{k=1}^{\infty}\sum_{\ell=1}^{\infty} \frac{d_{jk}b_\ell \Delta_{\ell k}}{\lambda_j + a_n}\phi_j(z).$$

Now $E(\Delta_{j_1 k_1}\Delta_{\ell_1 m_1}\Delta_{j_2 k_2}\Delta_{\ell_2 m_2}) = O(n^{-2})$, uniformly in the indicated indices; $\sum_\ell |b_\ell| < \infty$, since A.3 implies that $\beta > 1$; and $\sum_{k \geq 1}|d_{jk}| = O(j^{-\alpha/2})$, again by A.3. Therefore,

$$(E\|T^+ B_{n12}\|^4)^{1/2} = \left[E\left\{\sum_{j=1}^{\infty}\frac{1}{(\lambda_j + a_n)^2}\left(\sum_{k=1}^{\infty}\sum_{\ell=1}^{\infty}d_{jk}b_\ell\Delta_{\ell k}\right)^2\right\}^2\right]^{1/2}$$

(6.12)
$$= O\left\{\frac{1}{n}\sum_{j=1}^{\infty}\frac{1}{(\lambda_j+a_n)^2}\left(\sum_{k=1}^{\infty}\sum_{\ell=1}^{\infty}|d_{jk}||b_\ell|\right)^2\right\}$$

$$= O\left\{\frac{1}{n}\sum_{j=1}^{\infty}\frac{j^{-\alpha}}{(\lambda_j+a_n)^2}\right\} = O(n^{-(2\beta-1)/(2\beta+\alpha)}).$$

In view of (6.1) and (6.6),

(6.13) $E\|T^+\Delta\|^8 \leq \|T^+\|^8 E\|\Delta\|^8 = O(a_n^{-8} E\|\Delta\|^8) = O(\delta^4/a_n^8) = O(1).$

By (6.11), (6.13) and the Cauchy–Schwarz inequality,

(6.14)
$$(E\|T^+\Delta T^+ B_{n21}\|^4)^{1/2} \leq (E\|T^+\Delta\|^8 E\|T^+ B_{n21}\|^8)^{1/4}$$
$$= O(n^{-(2\beta-1)/(2\beta+\alpha)}).$$

Define

$$I_n(w) = \int\{\hat{f}_{XW}(x,w) - f_{XW}(x,w)\}g(x)\,dx,$$

$$J_n = \iint\{T^+(\hat{f}_{XW} - f_{XW})(z,w)\}^2\,dw\,dz.$$

Moment calculations show that $E\|I_n\|^8 = O(\delta^4)$ and $E(J_n^4) = O(\delta^4/a_n^8)$, and so by the Cauchy–Schwarz inequality,

(6.15)
$$(E\|T^+ B_{n3}\|^4)^{1/2} \leq \{E(\|I_n\|^4 J_n^2)\}^{1/2} \leq (E\|I_n\|^8 E J_n^4)^{1/4}$$
$$= O(\delta^2/a_n^2) = O(n^{-(2\beta-1)/(2\beta+\alpha)}).$$

It follows from (6.10)–(6.12), (6.14) and (6.15) that

(6.16) $$E\|\tilde{A}_{n31}\|^2 = O(n^{-(2\beta-1)/(2\beta+\alpha)}).$$

Now consider

(6.17)
$$\begin{aligned}(I+T^+\Delta)^{-1}T^+\Delta T^+ B_{n22} \\ = I(\|T^+\Delta\|\le\tfrac{1}{2})(I+T^+\Delta)^{-1}T^+\Delta T^+ B_{n22} \\ + I(\|T^+\Delta\|>\tfrac{1}{2})(I+T^+\Delta)^{-1}T^+\Delta T^+ B_{n22} \\ = H_{n1}+H_{n2},\end{aligned}$$

say. We first investigate $H_{n1}$.

If $\|T^+\Delta\|\le\tfrac{1}{2}$, then for some constant $D$ not depending on $\psi$, $\|(I+T^+\Delta)^{-1}\psi\|^2 \le D\|\psi\|^2$. Therefore, $\|H_{n1}\|^2 \le D\|T^+\Delta T^+ B_{n22}\|^2$. Some algebra shows that

$$T^+\Delta T^+ B_{n22}(z) = R_{n1}(z) + R_{n2}(z) + R_{n3}(z),$$

where

$$R_{n1}(z) = \int t^+(z,u)\{\hat f_{XW}(x,w_1)-f_{XW}(x,w_1)\}f_{XW}(u,w_1)$$
$$\times t^+(x,v)\{\hat f_{XW}(v,w_2)-E\hat f_{XW}(v,w_2)\}H(w_2)\,du\,dv\,dx\,dw_1\,dw_2,$$

$$R_{n2}(z) = \int t^+(z,u)\{\hat f_{XW}(u,w_1)-f_{XW}(u,w_1)\}f_{XW}(x,w_1)$$
$$\times t^+(x,v)\{\hat f_{XW}(v,w_2)-E\hat f_{XW}(v,w_2)\}H(w_2)\,du\,dv\,dx\,dw_1\,dw_2,$$

$$R_{n3}(z) = \int t^+(z,u)\{\hat f_{XW}(u,w_1)-f_{XW}(u,w_1)\}\{\hat f_{XW}(x,w_1)-f_{XW}(x,w_1)\}$$
$$\times t^+(x,v)\{\hat f_{XW}(v,w_2)-E\hat f_{XW}(v,w_2)\}H(w_2)\,du\,dv\,dx\,dw_1\,dw_2.$$

First we treat $R_{n1}$. Write $R_{n1}=R_{n11}+R_{n12}$, where

$$R_{n11}(z) = \int t^+(z,u)\{\hat f_{XW}(x,w_1)-E\hat f_{XW}(x,w_1)\}f_{XW}(u,w_1)$$
$$\times t^+(x,v)\{\hat f_{XW}(v,w_2)-E\hat f_{XW}(v,w_2)\}H(w_2)\,du\,dv\,dx\,dw_1\,dw_2,$$

$$R_{n12}(z) = \int t^+(z,u)\{E\hat f_{XW}(x,w_1)-f_{XW}(x,w_1)\}f_{XW}(u,w_1)$$
$$\times t^+(x,v)\{\hat f_{XW}(v,w_2)-E\hat f_{XW}(v,w_2)\}H(w_2)\,du\,dv\,dx\,dw_1\,dw_2.$$

By the Cauchy–Schwarz inequality,

$$\|R_{n11}\|^2 \le \int\left[\int t^+(z,u)\{\hat f_{XW}(x,w_1)-E\hat f_{XW}(x,w_1)\}\right.$$
$$\left.\times f_{XW}(u,w_1)\,du\,dw_1\right]^2 dx\,dz$$

$$\times \int [t^+(x,v)\{\hat{f}_{XW}(v,w_2) - E\hat{f}_{XW}(v,w_2)\}H(w_2)\,dv\,dw_2]^2\,dx$$

$$\equiv A_{n1}A_{n2},$$

say. Further application of the Cauchy–Schwarz inequality gives

(6.18) $$E\|R_{n11}\|^2 \leq \{(EA_{n1}^2)(EA_{n2}^2)\}^{1/2}.$$

Also,

(6.19) $$(EA_{n2}^2)^{1/2} = O\{(nha_n^2)^{-1}\}.$$

Now define $\delta_k(x) = \int \{\hat{f}_{XW}(x,w) - E\hat{f}_{XW}(x,w)\}\phi_k(w)\,dw$. Then

$$A_{n1} = \sum_{j=1}^{\infty}\sum_{k=1}^{\infty}\sum_{\ell=1}^{\infty} \frac{d_{jk}d_{j\ell}}{(\lambda_j + a_n)^2} \int \delta_k \delta_\ell,$$

from which it follows that

$$E(A_{n1}^2) = O\left\{\frac{1}{nh}\sum_{j=1}^{\infty}\sum_{k=1}^{\infty}\sum_{\ell=1}^{\infty} \frac{|d_{jk}||d_{j\ell}|}{(\lambda_j + a_n)^2}\right\} = O\left\{\frac{1}{nh}\sum_{j=1}^{\infty} \frac{j^{-\alpha}}{(\lambda_j + a_n)^2}\right\}$$

$$= O(h^{-1}n^{-(2\beta-1)/(2\beta+\alpha)}).$$

Combining this result with (6.18) and (6.19), we obtain

(6.20) $$E\|R_{n11}\|^2 = O\left(\frac{1}{nh^2 a_n^2}n^{-(2\beta-1)/(2\beta+\alpha)}\right) = O(n^{-(2\beta-1)/(2\beta+\alpha)}).$$

Calculations in the case of $R_{n12}$ are similar, as follows. We re-define

$$\delta_k(x) = \int \{E\hat{f}_{XW}(x,w) - f_{XW}(x,w)\}\phi_k(w)\,dw = O(h^r).$$

Therefore,

$$E\|R_{n12}\|^2 = O\left(\frac{h^{2r-1}}{a_n^2}n^{-(2\beta-1)/(2\beta+\alpha)}\right) = O(n^{-(2\beta-1)/(2\beta+\alpha)}).$$

Combining this result with (6.20), we deduce that

(6.21) $$E\|R_{n1}\|^2 = O(n^{-(2\beta-1)/(2\beta+\alpha)}).$$

Next we treat $R_{n2}$. Re-define $A_{n1}$ and $A_{n2}$ by

$$R_{n2}(z)^2 \leq \int \left[\int t^+(z,u)\{\hat{f}_{XW}(u,w_1) - f_{XW}(u,w_1)\}\,du\right]^2 dw_1$$

$$\times \int \left[\int f_{XW}(x,w_1)t^+(x,v)\right.$$

$$\left.\times \{\hat{f}_{XW}(v,w_2) - E\hat{f}_{XW}(v,w_2)\}H(w_2)\,dv\,dx\,dw_2\right]^2 dw_1$$

$$= A_{n1}(z)A_{n2}(z).$$

Furthermore,

$$(6.22) \qquad (E\|A_{n1}\|^2)^{1/2} = O\left(\frac{h^{2r}}{a_n^2} + \frac{1}{nha_n^2}\right).$$

Defining $\delta_{jk} = \int\{\hat{f}_{XW}(x,w) - E\hat{f}(_{XW}(x,w)\}\phi_j(x)\phi_k(w)\,dx\,dw$ and $h_j = \int H\phi_j$, we have

$$\int A_{n2} = \sum_{k=1}^{\infty}\left(\sum_{j=1}^{\infty}\sum_{\ell=1}^{\infty}\frac{d_{jk}\delta_{j\ell}h_\ell}{\lambda_j + a_n}\right)^2 = \sum_{j=1}^{\infty}\sum_{\ell=1}^{\infty}\sum_{s=1}^{\infty}\frac{\lambda_j \delta_{j\ell}\delta_{js}h_\ell h_s}{(\lambda_r + a_n)^2}.$$

Therefore,

$$(E\|A_{n2}\|^2)^{1/2} = O\left\{n^{-1}\sum_{j=1}^{\infty}\frac{\lambda_j}{(\lambda_j+a_n)^2}\right\} = O(n^{-(2\beta-1)/(2\beta+\alpha)}).$$

This result and (6.22) give

$$(6.23) \qquad E\|R_{n2}\|^2 = O(n^{-(2\beta-1)/(2\beta+\alpha)}).$$

Next we treat $R_{n3}$. Note that $R_{n3}(z)^2 \leq A_{n1}(z)A_{n2}(z)$, where we re-define

$$A_{n1}(z) = \int\left[\int t^+(z,u)\{\hat{f}_{XW}(u,w_1) - f_{XW}(u,w_1)\}\right.$$
$$\left. \times \{\hat{f}_{XW}(x,w_1) - f_{XW}(x,w_1)\}\,du\,dw_1\right]^2 dx,$$

$$A_{n2}(z) = \int\left[\int t^+(x,v)\{\hat{f}_{XW}(v,w_2) - E\hat{f}_{XW}(v,w_2)\}H(w_2)\,dv\,dw_2\right]^2 dx.$$

Therefore,

$$E\|R_{n3}\|^2 = O\left(\frac{1}{n^3h^5a_n^4} + \frac{h^{4r}}{nha_n^2}\right) = O(n^{-(2\beta-1)/(2\beta+\alpha)}).$$

Combining this result with (6.21) and (6.23), and recalling the definition of $H_{n1}$ at (6.17), we deduce that

$$(6.24) \qquad E\|H_{n1}\|^2 = O(n^{-(2\beta-1)/(2\beta+\alpha)}).$$

Now we consider $H_{n2}$. We have

$$\|(I + T^+\Delta)^{-1}\psi\| = \|\widehat{T}^+(T + a_n I)\psi\|$$
$$\leq \|\widehat{T}^+\|\|T + a_n I\|\|\psi\|$$
$$\leq \text{const.}\, a_n I^{-1}\|\psi\|.$$

Therefore,

$$\|H_{n2}\|^2 \leq \text{const.}\, a_n^{-2} I(\|T^+\Delta\| > \tfrac{1}{2})\|T^+\Delta T^+ B_{n22}\|^2,$$

and so by the Cauchy–Schwarz inequality,

$$E\|H_{n2}\|^2 \leq \text{const.}\, a_n^{-2} P(\|T^+\Delta\| > \tfrac{1}{2})^{1/2} (E\|T^+\Delta T^+ B_{n22}\|^4)^{1/2}.$$

We shall prove shortly that, for all $\ell > 0$,

(6.25) $$P(\|T^+\Delta\| > \tfrac{1}{2}) = O\{(\delta/a_n^2)^\ell\}.$$

Moreover,

(6.26) $$E\|T^+ B_{n22}\|^8 \leq \|T^+\|^8 E\|B_{n22}\|^8 = O(a_n^{-8} E\|B_{n22}\|^8)$$
$$= O\left\{a_n^{-8}\left(\int B_{n22}^2\right)^4\right\} = O\{(nha_n^2)^{-4}\},$$

the last identity following by moment calculations similar to those leading to (6.6). Combining (6.13) and (6.26), and applying the Cauchy–Schwarz inequality, we deduce that

$$(E\|T^+\Delta T^+ B_{n22}\|^4)^{1/2} \leq (E\|T^+\Delta\|^8 E\|T^+ B_{n22}\|^8)^{1/4} = O\{(\delta/a_n^2)(nha_n^2)^{-1}\}.$$

Using this result together with (6.25), and choosing $\ell$ sufficiently large, we obtain

$$E\|H_{n2}\|^2 = O\{(\delta/a_n^2)^{1+(\ell/2)}(nha_n^2)^{-1}\} = O(n^{-(2\beta-1)/(2\beta+\alpha)}).$$

Combining this result with (6.17) and (6.24), we obtain

$$E\|(I+T^+\Delta)^{-1}T^+\Delta T^+ B_{n22}\|^2 = O(E\|H_{n1}\|^2 + E\|H_{n2}\|^2)$$
$$= O(n^{-(2\beta-1)/(2\beta+\alpha)}).$$

Result (6.3) for $j = 3$ follows from this formula and (6.16).

Next we derive (6.3) for $j = 4$. Since $\widehat{T}^+ - T^+ = -(I+T^+\Delta)^{-1}T^+\Delta T^+$ and $I - \widehat{T}^+ T = -(I+T^+\Delta)^{-1}T^+\Delta$, then

$$A_{n4} = -(I+T^+\Delta)^{-1}T^+\Delta(A_{n2} - T^+ B_{n2}).$$

The arguments leading to (6.3) with $j = 2$, and (6.15), may be used to prove that

$$\eta \equiv \{(\delta^2/a_n)^4 E\|A_{n2}\|^4 + E\|T^+\Delta T^+ B_{n2}\|^4\}^{1/2} = O(n^{-(2\beta-1)/(2\beta+\alpha)}).$$

Therefore, by (6.7), (6.8) and the Cauchy–Schwarz inequality,

$$E\|A_{n4}\|^2 \leq 2\,\{E\|(I+T^+\Delta)^{-1}T^+\Delta\|^4 E\|A_{n2}\|^4\}^{1/2}$$
$$+ 2\{E\|(I+T^+\Delta)^{-1}\|^4 E\|T^+\Delta T^+ B_{n2}\|^4\}^{1/2}$$
$$= O(\eta) = O(n^{-(2\beta-1)/(2\beta+\alpha)}).$$

This proves (6.3) for $j = 4$.

It remains to derive (6.7). Let $\psi \in L_2[0,1]$. Then, for constants not depending on $\psi$, if $\|T^+\Delta\| \leq \frac{1}{2}$,

$$\|(I + T^+\Delta)^{-1}\psi\| \leq \text{const.}\,\|\psi\|$$

and, without any constraint on $\|T^+\Delta\|$,

$$\begin{aligned}\|(I + T^+\Delta)^{-1}\psi\| &= \|\widehat{T}^+(T + a_n I)\psi\| \\ &\leq \|\widehat{T}^+\|\|T + a_n I\|\|\psi\| \\ &\leq \text{const.}\, a_n^{-1}\|\psi\|.\end{aligned}$$

Therefore,

$$\|(I + T^+\Delta)^{-1}\| \leq \text{const.}\{1 + a_n^{-1} I(\|T^+\Delta\| > \tfrac{1}{2})\}.$$

Hence, noting (6.6), and employing Markov's inequality to bound $P(\|T^+\Delta\| > \frac{1}{2})$, we deduce that, for each fixed $k, \ell > 0$,

(6.27)
$$\begin{aligned}E\{\|(I + T^+\Delta)^{-1}\|^k\} &\leq \text{const.}\{1 + a_n^{-k} P(\|T^+\Delta\| > \tfrac{1}{2})\} \\ &\leq \text{const.}\{1 + a_n^{-k} E(\|T^+\Delta\|^{2\ell})\} \\ &\leq \text{const.}\{1 + a_n^{-k-2\ell} E(\|\Delta\|^{2\ell})\} \\ &\leq \text{const.}(1 + a_n^{-k-2\ell}\delta^\ell) \\ &= \text{const.}\{1 + a_n^{-k}(\delta/a_n^2)^\ell\},\end{aligned}$$

where the constants depend on $k$ and $\ell$ but not on $n$. If $k$ is given, then we may choose $\ell = \ell(k)$ so large that $a_n^{-k}(\delta/a_n^2)^\ell \to 0$ as $n \to \infty$, and so (6.7) follows from (6.27). This argument also gives (6.25).

6.2. *Proof of Theorem* 4.2. Put $\bar{p} = (p_1, \ldots, p_m)'$, where $p_j = E_G\{g(\widetilde{X}) \times \chi_j(\widetilde{W})\} = E_G\{Y\chi_j(\widetilde{W})\}$. Let $\gamma = (\gamma_j)$ and $p = (p_j)$ denote infinite column vectors, and let $\bar{Q}$ be the $m \times m$ upper left-hand sub-matrix of $Q$. Since $p = Q\gamma$, then $p_j = p_j(G) = O(j^{-(2\beta+\alpha)/2})$, uniformly in $G \in \mathcal{H}$, as $j \to \infty$. Therefore, $(\bar{Q}'p)_i = O(i^{-(\alpha+\beta)})$, uniformly in $1 \leq i \leq m$, $n \geq 1$ and $G \in \mathcal{H}$. This result will be used below without further reference.

Put $\bar{M} = \bar{Q}\bar{Q}' + a_n I_m$ and $\widehat{M} = \widehat{Q}\widehat{Q}' + a_n I_m$. It may be deduced from the definition of $\mathcal{H}$ that the bounds on $|q_{jk}|$ and $|q_{jk}^{(-1)}|$ in that definition apply too to the $(j,k)$th elements of $\bar{M}$ and $\bar{M}^{-1}$, respectively, provided we replace $\alpha$ by $2\alpha$ and alter the constants $C_1$ and $C_2$ (retaining their positivity, of course). The bounds are valid uniformly in $1 \leq j, k \leq m$ and $n \geq 1$, and permit it to be proved that

$$(\bar{M}^{-1}\bar{Q}'\bar{p})_j = \{(Q'Q)^{-1}Q'p\}_j + O(m^{-\beta}) = \gamma_j + O(m^{-\beta}),$$

uniformly in $1 \leq i \leq m$, $n \geq 1$ and distributions of $G \in \mathcal{H}$. Note too that
$$\widehat{M}^{-1} = \{I + \bar{M}^{-1}(\widehat{M} - \bar{M})\}^{-1}\bar{M}^{-1},$$
$$\widehat{M}^{-1}\widehat{Q}'\hat{p} - \bar{M}^{-1}\bar{Q}'\bar{p} = \{\bar{M}^{-1} + (\widehat{M}^{-1} - \bar{M}^{-1})\}(\widehat{Q}'\hat{p} - \bar{Q}'\bar{p})$$
$$+ (\widehat{M}^{-1} - \bar{M}^{-1})\bar{Q}'\bar{p}.$$

From these properties it may be shown that

(6.28)
$$E_G\left\{\sum_{j=1}^m (\tilde{\gamma}_j - \gamma_j)^2\right\}$$
$$= O\left\{E_G\left(\sum_{i=1}^m [\{\bar{M}^{-1}(\widehat{Q}'\hat{p} - \bar{Q}'\bar{p})\}_j]^2\right)\right.$$
$$+ E_G\left(\sum_{i=1}^m [\{\bar{M}^{-1}(\widehat{M} - \bar{M})\bar{M}^{-1}(\widehat{Q}'\hat{p} - \bar{Q}'\bar{p})\}_j]^2\right)$$
$$\left. + E_G\left(\sum_{i=1}^m [\{\bar{M}^{-1}(\widehat{M} - \bar{M})\bar{M}^{-1}\bar{Q}'\bar{p}\}_j]^2\right) + m^{1-2\beta}\right\},$$

uniformly in $G \in \mathcal{H}$.

It may be proved by Taylor expansion arguments, involving approximating $\widehat{W}_i = \widehat{F}_W(W_i)$ by $\widetilde{W}_i = F_W(W_i)$, and analogously for $\widehat{X}_i$ and $\widetilde{X}_i$, that, for each $r, \varepsilon > 0$,

(6.29)
$$\max_{1 \leq j,k \leq n^{(1/2)-\varepsilon}} \sup_{G \in \mathcal{H}} E_G|\hat{q}_{jk} - q_{jk}|^r = O(n^{-r/2}),$$

(6.30)
$$\max_{1 \leq j \leq n^{(1/2)-\varepsilon}} \sup_{G \in \mathcal{H}} E_G(\hat{p}_j - p_j)^2 = O(n^{-1}).$$

Rather standard, but tedious, moment calculations, using (6.29) and (6.30), may be employed to show that each of the expected values on the right-hand side of (6.28) equals $O(n^{-1}m^{\alpha+1})$, uniformly in $G \in \mathcal{H}$. Therefore,

(6.31)
$$\sup_{G \in \mathcal{H}} \sum_{j=1}^m E_G\{(\tilde{\gamma}_j - \gamma_j)^2\} = O(n^{-1}m^{\alpha+1} + m^{1-2\beta})$$
$$= O(n^{-(2\beta-1)/(2\beta+\alpha)}).$$

It follows from the definition of $\mathcal{H}$ that $\sum_{j>m} \gamma_j^2 = O(m^{1-2\beta})$, uniformly in $G \in \mathcal{H}$. This result and (6.31) imply that

$$\int E_G(\bar{g} - g)^2 = \sum_{j=1}^m E_G(\tilde{\gamma}_j - \gamma_j)^2 + \sum_{j=m+1}^\infty \gamma_j^2 = O(n^{-(2\beta-1)/(2\beta+\alpha)})$$

uniformly in $G \in \mathcal{H}$, completing the proof of the theorem.

6.3. *Proof of Theorem* 4.4. For simplicity, we deal only with the orthogonal series setting, discussed in Section 4.2. We may assume the following: $\phi_j \equiv \chi_j$, $\phi_1 \equiv 1$ and $\phi_{j+1}(x) = 2^{-1/2}\cos(j\pi x)$, for $j \geq 1$; the marginal distributions of $X$ and $W$ are uniform on the unit interval; and

$$
\begin{aligned}
f_{XW}(x,w) &= \sum_{j=1}^{\infty} j^{-\alpha/2}\phi_j(x)\phi_j(w), \\
Y &= \sum_{j=m+1}^{2m} \theta_j j^{-(2\beta+\alpha)/2}\phi_j(W) + V,
\end{aligned}
\tag{6.32}
$$

where $m$ equals the integer part of $n^{1/(2\beta+\alpha)}$, the $\theta_j$'s are all either 0 or 1, and $V$ is Normal $N(0,1)$, independent of $(X,W)$.

The function $g$ implied by (6.32) is $g(x) = \sum_{m+1 \leq j \leq 2m} \theta_j j^{-\beta}\phi_j(x)$. Note too that if $\tilde{g}$ is an estimator of $g$, then

$$
\tilde{\theta}_j = j^{\beta}\int \tilde{g}\phi_j \tag{6.33}
$$

may be viewed as an estimator of $\theta_j$.

A standard argument based on the Neyman–Pearson lemma shows that

$$
\liminf_{n\to\infty} \inf_{m+1\leq j\leq 2m} \inf_{\check{\theta}_j} \sup{}^* E(\check{\theta}_j - \theta_j)^2 > 0, \tag{6.34}
$$

where $\sup^*$ denotes the supremum over all $2^m$ different distributions of $(X,W,Y)$ obtained by taking different choices of $\theta_{m+1}, \ldots, \theta_{2m}$ in (6.32), and $\inf_{\check{\theta}_j}$ represents the infimum over all measurable functions $\check{\theta}_j$ of the data. To derive (6.34), it suffices to take $\check{\theta}_j$ to be the likelihood-ratio rule for distinguishing between $\theta_j = 0$ and $\theta_j = 1$, and work through a little asymptotic theory to obtain the version of (6.34) when "$\inf_{\check{\theta}_j}$" is omitted from the left-hand side.

Therefore, if $\tilde{g}$ is given, and $\tilde{\theta}_{m+1}, \ldots, \tilde{\theta}_{2m}$ are the estimators of $\theta_{m+1}, \ldots, \theta_{2m}$, respectively, derived from $\tilde{g}$ as suggested at (6.33), then

$$
\begin{aligned}
\sup{}^* \int (\tilde{g} - g)^2 &= \sup{}^* \sum_{j=m+1}^{2m} E(\tilde{\theta}_j - \theta_j)^2 j^{-2\beta} \\
&\geq \text{const.} \sum_{j=m+1}^{2m} j^{-2\beta} \\
&\geq \text{const.} \, j^{-(2\beta-1)/(2\beta+\alpha)},
\end{aligned}
$$

where the constants do not depend on choice of $\tilde{g}$. This proves the theorem.

Centre for Mathematics
  and Its Applications
Australian National University
Canberra, ACT 0200
Australia
e-mail: halpstat@maths.anu.edu.au

Department of Economics
Anderson Hall
Northwestern University
2001 Sheridan Road
Evanston, Illinois 60208-2600
USA